\documentclass[11pt,a4]{amsart}
\makeatletter
\renewcommand{\section}{\@startsection{section}{1}{0pt}{20pt}{6pt}{\large\bfseries}}
\makeatother

\usepackage{enumerate, pdfsync}
\usepackage{color}
\usepackage{xcolor}
\usepackage{graphicx}

\numberwithin{equation}{section}

\theoremstyle{plain}
  \newtheorem{thm}{Theorem}[section]
  
  \newtheorem{lemma}[thm]{Lemma}

\theoremstyle{definition}
  
\newtheorem{exa}[thm]{Example}

\newcommand{\R}{\mathbb{R}}

\newcommand{\C}{\mathbb{C}}
\newcommand{\Q}{{\rm  P}}
\newcommand{\E}{\mathbb{E}}

\renewcommand{\P}{{\mathbb{P}}}

\newcommand{\iex}{A}

\newcommand{\I}{\mathcal{I}}
\newcommand{\Ip}{\mathcal{I}_{\psi,\alpha}}

\newcommand{\ant}{a_n(\psi;\alpha)}

\newcommand{\Id}[1]{{{\mathbb{I}}}_{\{#1\}}}

\newcommand{\Oc}[2]{O_q^{\mathcal{T}_\alpha\psi}(#1;#2)}
\begin{document}

\bibliographystyle{plain}

\title{A Ciesielski-Taylor type identity for positive self-similar Markov processes}
\author{A. E. Kyprianou}
\address{Department of Mathematical Sciences
University of Bath, Bath BA2 7AY, UK} \email{a.kyprianou@bath.ac.uk}

\author{P. Patie}
   \address{D\'epartement de Math\'ematiques, Universit\'e Libre de Bruxelles\\
Boulevard du Triomphe,  B-1050 Bruxelles.}
\email{ppatie@ac.ulb.be}

\begin{abstract}
The aim of this note
is to give a straightforward proof of a general version of the Ciesielski-Taylor identity for positive self-similar Markov processes of the spectrally negative type which umbrellas all previously known Ciesielski-Taylor identities within the latter class. The approach makes use of three fundamental features. Firstly a new transformation which maps a subset of the family of Laplace exponents of spectrally negative L\'evy processes into itself. Secondly some classical features of fluctuation theory for spectrally negative L\'evy processes (see eg.~\cite{Kyprianou-Palmowski-05}) as well as more recent fluctuation identities for positive self-similar Markov processes found in  Patie \cite{Patie-06c}.

\bigskip

\noindent \textsc{R\'esum\'e.} L'objectif principal de ce papier est de donner une preuve d'une version g\'en\'erale de l'identit\'e de Ciesielski-Taylor pour la famille de processus positifs auto-similaires markoviens et de type spectralement n\'egatif, ce qui nous permet d'unifier l'ensemble des r\'esultats d\'ej\`a connus pour cette famille. Notre preuve s'appuie sur trois concepts importants. Tout d'abord, nous introduisons une famille de transformations qui associe l'ensemble des exposants de Laplace de processus de L\'evy spectralement n\'egatifs \`a lui-m\^eme. Ensuite nous combinons des r\'esultats emprunt\'es \`a   la th\'eorie des fluctuations des processus de L\'evy  spectralement n\'egatifs (voir e.g.~\cite{Kyprianou-Palmowski-05}) et \`a celles  des processus  positifs auto-similaires markoviens  spectralement n\'egatif
\'elabor\'ees plus r\'ecemment par Patie \cite{Patie-06c}.

\bigskip


\bigskip

\end{abstract}

\date{}

\maketitle


\section{Introduction}
Suppose that $(X, Q^{(\nu)})$  is a Bessel process $X$ starting from $0$ with dimension $\nu>0$. That is to say,  the $[0,\infty)$-valued diffusion whose infinitesimal generator is given by
\[
L_\nu f(x) = \frac{1}{2}f''(x) +\frac{\nu-1}{2x}f'(x)
\]
on $(0,\infty)$ for $f\in C^2(0,\infty)$ with instantaneous reflection at $0$ when $\nu\in (0,2)$ (i.e. $f'(0^+)=0$) and when $\nu\geq 2$ the origin is  an entrance-non-exit boundary point.
For these processes, Ciesielski and Taylor \cite{Ciesielski-Taylor-62} observed that the following curious identity in distribution holds. For  $a>0$ and {\it integer} $\nu>0$,
\begin{equation} \label{eq:ict}
\left(T_a, Q^{(\nu)}\right) \stackrel{(d)}{=}
\left(\int_0^{\infty}\Id{X_s\leq a} ds, Q^{(\nu+2)}\right)
\end{equation}
where \[
 T_a
=\inf\{s \geq 0; \: X_s = a\}.\]  They proved this relationship by
showing that the densities of both random variables coincide. Getoor
and Sharpe \cite{Getoor-Sharpe-79} extended this identity to any dimension $\nu>0$ by means of the Laplace
transform and recurrence relationships for Bessel functions.  Biane
\cite{Biane-85} has generalized this identity in law to one dimensional diffusions
by appealing to the Feynman-Kac formula for the Laplace transforms of the path functionals involved and an analytical  manipulation of the associated  infinitesimal generators. 
Yor
\cite{Yor-91} offered a probabilistic explanation by using the
occupation times formula and Ray-Knight theorems. Finally, Carmona
et al.~\cite[Theorem 4.8]{Carmona-Petit-Yor-98} proved a similar
identity, in terms of the confluent hypergeometric function,
 for a self-similar `saw tooth' process. There is also a  Ciesielski-Taylor type identity for  spectrally negative L\'evy processes which is to be found in a short remark of  Bertoin \cite{BertoinCT}.

 In the majority of the aforementioned cases, the underlying stochastic processes are examples of positive self-similar Markov processes of the spectrally negative type.  Recall that Lamperti \cite{Lamperti-72} showed that, for any $x\in\R$, there
exists a one to one mapping between $\P_x$, the law of a generic L\'evy process (possibly killed at an independent and exponentially distributed time), say $\xi = (\xi_t : t\geq 0)$, starting from $x$, and the law
$\Q_{e^{x}}$ of an $\alpha$-self-similar positive Markov process, say  $X = (X_t : t\geq 0)$, starting from $e^{x}$ and killed on first hitting zero. The latter process is a $[0,\infty)$-valued Feller process which enjoys the following
$\alpha$-self-similarity property, for any $\alpha,x> 0,$ and $c>0$,
\begin{equation} \label{eq:self}
\left((X_{t})_{t\geq0}, \Q_{cx}\right)
\stackrel{(d)}{=}\left((cX_{c^{-\alpha} t})_{t\geq0}, \Q_{x}\right).
\end{equation}
Specifically, Lamperti proved that $X$ can be constructed from
$\xi$ via the relation
\begin{equation} \label{eq:lamp_transf}
\log X_t  =  \xi_{\iex_{t }}, \: 0\leq t< \zeta,
\end{equation}
 where $\zeta = \inf\{t>0 : X_t =0\}$ and
\begin{equation}
 \iex_t = \inf
\{ s \geq 0; \:  \int_0^s e^{\alpha \xi_u} \: du
> t \}.
\label{A_t}
\end{equation}
In this paper we are predominantly interested in the case that
 $\xi$ is a (possibly killed) spectrally negative L\'evy process; that is to say,   $X$ is a positive self-similar Markov process of the spectrally negative type. For this class of driving spectrally negative L\'evy processes it is  known that  $\mathbb{E}(\xi_1)\in[-\infty, \infty)$ and when there is no killing and $\mathbb{E}(\xi_1)\geq 0$
one may extend the definition of $X$ to include the case that it is issued from the origin by establishing its entrance law ${\rm  P}_0$ as the weak limit with respect to the Skorohod topology of ${\rm  P}_{x}$ as $x\downarrow 0$; see Bertoin and Yor  \cite{Bertoin-Yor-02-b} and  Chaumont and Caballero \cite{Ca-Ch}.
 We also recall that when $\mathbb{E}(\xi_1)<0$ (resp. $\xi$ is killed) then the boundary state $0$ is reached continuously (resp.~by a jump). In these two cases, one cannot construct an entrance law, however, Rivero \cite{Rivero1} and Fitzsimmons \cite{Fitz}, show that it is possible instead to construct a unique recurrent extension on $[0,\infty)$ such that paths leave $0$ continuously,  thereby giving  a meaning to ${\rm  P}_0$, if and only if  there exists a $\theta\in(0,\alpha)$ such that $\mathbb{E}(e^{\theta\xi_1}) = 1$.
Moreover, for these cases, the recurrent extension on $[0,\infty)$ is unique.


The object of this paper is to establish a new general Ciesielski-Taylor type identity for  the aforementioned class of  self-similar Markov processes of the spectrally negative type issued from the origin. Our identity will umbrella  all of the  known examples within this class. The basis of our new identity will be the blend of a new transformation which maps a subset of the family of Laplace exponents of spectrally negative L\'evy processes into itself together with  some classical features of fluctuation theory for spectrally negative L\'evy processes as well as more recent fluctuation identities for positive self-similar Markov processes. Although we appeal to the principle of matching Laplace transforms in order to obtain our distributional identity, we consider our proof to be largely probabilistic and quite straighforward in its nature. We  make predominant  use of spectral negativity  and the strong Markov property, avoiding the use of the Feynman kac formula and subsequent integro-differential equations that follow thereof.

The remaining part of the paper is organized as follows. In the next section, we first introduce preliminary notation as well as the family of transformations $\mathcal{T}_\beta$ acting on Laplace exponents of spectrally negative L\'evy processes. In the section thereafter we state and prove our new Ciesielski-Taylor type identity. Finally we conclude with some examples including some discussion on how our technique relates to possible alternative proofs which make use of the Feynman-Kac formula.


\section{The transformation $\mathcal{T}_\beta$}

Recall that a killed L\'evy process should be understood as the process which executes the path of a L\'evy process up to an independent and exponentially distributed random time at which point it is sent to a cemetery state which is taken to be $+\infty$. {\it Henceforth, when referring to a L\'evy process, we shall implicitly understand that the possibility of killing is allowed.}
For any spectrally negative L\'evy process, henceforth denoted by  $\xi = (\xi_t, t\geq 0)$, whenever it exists we define the Laplace exponent  by
\[
\psi(u) = \log \mathbb{E}_0(\exp\{u \xi_1\}).
\]
It is a well established fact that the latter Laplace exponent is well defined and strictly convex on $[0,\infty)$, see for example Bertoin \cite{Bertoin-96}. 

  Next we  introduce the family of  transformations $\mathcal{T}=(\mathcal{T}_\beta)_{\beta\geq 0}$  acting on  Laplace exponents of spectrally negative L\'evy processes.

\begin{lemma}\label{lem:transf}
For each fixed  $\beta\geq 0$, define the linear transformation
\[
\mathcal{T}_\beta \psi(u) = \frac{u}{u+\beta}\psi(u+\beta), \, u \geq 0,
\]
for all  Laplace exponents $\psi$ of a spectrally negative L\'evy processes. 
Then $\mathcal{T}_\beta\psi$ is the Laplace exponent of a spectrally negative L\'evy process which is without killing whenever $\beta> 0$. Moreover,   the operator $\mathcal{T}$ satisfies  the composition property $\mathcal{T}_\beta \circ \mathcal{T}_\gamma = \mathcal{T}_{\beta+\gamma}$, for any $\beta,\gamma\geq 0$.
\end{lemma}


\begin{proof}[Proof of Lemma \ref{lem:transf}] First note that it is trivial that $\mathcal{T}_\beta$
is a linear transformation. Suppose first that $\psi(u) =\psi^*(u) - q$ where $q>0$ and $\psi^*(u) $ is the Laplace exponent of a spectrally negative L\'evy process with no killing. Then for $u,\beta\geq 0$,
\begin{equation}
\mathcal{T}_\beta \psi(u) =\mathcal{T}_\beta\psi^*(u) -q\frac{u}{u+\beta} =\mathcal{T}_\beta\psi^*(u) -q\int_0^\infty (1-e^{-u x})\beta e^{-\beta x}dx,
\label{q}
\end{equation}
showing that the effect of $\mathcal{T}_\beta$ on killed spectrally negative L\'evy processes is to subtract an additional compound Poisson subordinator with exponentially distributed jumps from the  transformed process without killing.
Hence
it suffices to prove the first claim for exponents which do not have a killing term.

To this end we assume henceforth that $q=0$ and we
 define the Esscher transformation on functions $f$ by  $\mathcal{E}_\beta f(u) = f(u+\beta) - f(\beta)$ whenever it makes sense. 
 A straightforward computation shows that
\begin{equation}
\mathcal{T}_\beta \psi (u)= \mathcal{E}_\beta \psi(u) - \beta \mathcal{E}_\beta \phi(u)
\label{T}
\end{equation}
where $\phi(u) = \psi(u)/u$.
It is well known that $ \mathcal{E}_\beta \psi(u) $ is the Laplace exponent of a spectrally negative L\'evy process with no killing and hence the proof is complete if we can show that $\mathcal{E}_\beta \phi(u)$ is  the Laplace exponent of a subordinator without killing. Indeed this would show that $\mathcal{T}_\beta \psi (u)$ corresponds to the independent sum of an Esscher transformed version of the original spectrally negative L\'evy process and the negative of a subordinator.

Thanks to the Wiener-Hopf factorization we may always write
$
\psi(u) = (u - \theta)\varphi(u)
$
where $\varphi$ is the Laplace exponent of the (possibly killed) subordinator which plays the role of  the descending ladder height process of the spectrally negative L\'evy process associated with $\psi$ and $\theta$ is the largest root in $[0,\infty)$ of the equation $\psi(\theta) = 0$.
If $\theta = 0$ then the proof is complete as $\phi = \varphi$ and hence $\mathcal{E}_\beta\phi(u)$ is the Esscher transform of a subordinator exponent which is again the Laplace exponent of a subordinator.

 Now assume that  $\theta>0$ and note that this is the case if and only if the aforementioned L\'evy process drifts $-\infty$ which implies that $\varphi$ has no killing component. (See Chapter 8 of Kyprianou \cite{Kyprianou-06}). Hence, using the same idea as in  (\ref{T}) again, we have
\begin{equation}\label{lastterm}
\phi(u)  = \frac{u-\theta}{u}\varphi(u)  = \mathcal{T}_\theta\varphi(u-\theta) =
 \mathcal{E}_\theta \varphi(u-\theta) - \theta \mathcal{E}_\theta \eta(u-\theta)
\end{equation}
where, if $\theta>0$,
\[
\eta(u) = \varphi(u)/u = {\rm d} + \int_0^\infty e^{- u x}\nu(x,\infty)dx,
\]
${\rm d\geq 0}$ and $\nu$ is a measure on $(0,\infty)$ which satisfies $\int_0^\infty(1\wedge x)\nu(dx)<\infty$. 
Thus
\[
\mathcal{E}_\beta \phi(u) = \mathcal{E}_{\beta + \theta} \varphi(u-\theta) - \theta \mathcal{E}_{\beta +\theta} \eta(u-\theta)
=
 \mathcal{E}_{\beta} \varphi(u)  + \theta\int_0^\infty (1 - e^{-u x})e^{-\beta x}\nu(x,\infty)dx
\]
which is indeed the Laplace exponent of a subordinator without killing thanks to the fact that the Esscher transform on a subordinator produces a subordinator without killing.

The final claim in the statement of the lemma is easily verified from the definition of the transformation.
\end{proof}

\section{Ciesielski-Taylor type identity}
Fix $\alpha>0$ and as in the previous section,  $\psi$ will denote the Laplace exponent of a given spectrally negative L\'evy process. As we wish to associate more clearly the underlying L\'evy process with each positive self-similar Markov process, we shall work with the modified notation $\mathbb{P}^\psi_\cdot$ and ${\rm P}^\psi_\cdot$ with the obvious choice of notation for their respective expectation operators. We emphasize here again for clarity, in the case that $x=0$, we understand ${\rm P}^\psi_0$ to be the law of the recurrent extension of $X$ when $\psi'(0^+)<0$ or $\psi$ has a killing term  and otherwise when $\psi'(0^+)\geq 0$ and $\psi$ has no killing term,  it is understood to be the  entrance law.

 Next we  introduce more notation taken from Patie \cite{Patie-06c}. Define  for non-negative integers $n$
\begin{equation*}
\ant^{-1}=\prod_{k=1}^n \psi(\alpha k) , \quad a_0=1,
\end{equation*}
and we introduce the entire function $\Ip$ which admits the series
representation
\begin{equation*}
\Ip(z)=\sum_{n=0}^{\infty} \ant z^{n}, \quad
 z
\in  \C.
\end{equation*}
It is important to note that whenever $\theta$, the largest root of the equation $\psi(\theta)=0$, satisfies $\theta<\alpha$, it follows that all of the coefficients in the definition of $\Ip(z)$ are strictly positive.



\begin{thm} \label{thm:ct} Fix $\alpha>0$. Suppose that $\psi$ is the Laplace exponent of a (possibly killed) spectrally negative L\'evy process. Assume that $\theta$,  the largest root in $[0,\infty)$ of the equation $\psi(\theta) =0$, satisfies $\theta<\alpha$. Then 
for any $a>0$, the following Ciesielski-Taylor type identity  in law
\begin{equation}
\left(T_a, \Q^{\psi}_0\right)
\stackrel{(d)}{=} \left(\int_0^{\infty}\Id{X_s\leq a} ds,
\Q^{\mathcal{T}_{\alpha}\psi}_0\right)
\end{equation}
holds.
Moreover, 
both random variables under their respective measures are self-decomposable and
their Laplace transforms in $q>0$ are both equal to
\[
\frac{1}{\mathcal{I}_{\psi,\alpha}(q a^\alpha)}.
\]

\end{thm}

The remaining part of this section is devoted to the proof of Theorem \ref{thm:ct}. Let us start by recalling that, in the case $\psi(0)=0$ and $\psi(\theta)=0$ with $\theta\in[0,\alpha)$, the Laplace transform of $$T_a:=\inf\{t>0: X_t =a\}$$
has been characterized by Patie \cite[Theorem 2.1]{Patie-06c} as follows. For   $0 \leq x \leq a$ and $q\geq0$, we have
\begin{equation}  \label{eq:lt}
 {\rm{E}}^{\psi}_x \left[e^{-q T_a } \right]
 = \frac{\Ip(qx^{\alpha})}{\Ip(qa^{\alpha})}.
\end{equation}
This identity will be used in the proof of Theorem \ref{thm:ct}, however it does not cover the case when $\psi(0)<0$ (i.e.~the underlying L\'evy process is killed) and $\psi(\theta)=0$ for $\theta<\alpha$. The following lemma fills this gap.

\begin{lemma}\label{lem:lt}
Suppose that $\psi(0)<0$ and the root of $\psi(\theta) = 0$ satisfies $\theta\in[0,\alpha)$.
 Then, for any $x>0$,
\[\lim_{t \downarrow 0} \frac{\Q_x(T_0\leq t)}{t} = -\psi(0)x^{-\alpha}.\]
Moreover, the infinitesimal generator $L^{\psi}$ of the recurrent extension has the following form
\[L^{\psi} f(x)= x^{-\alpha} L^{\xi} (f\circ \exp)(\log(x))-x^{-\alpha}\psi(0)f(0),\quad x>0,\]
for at least functions $f$ such that $f(x),xf'(x),x^2f''(x)$ are continuous on $[0,\infty)$ with $\lim_{x\downarrow 0}x^{\theta-1}f'(x)=0$ and $L^{\xi}$ is the infinitesimal generator of the killed L\'evy process $\xi$.
Finally, for any $0 \leq x \leq a$ and  $q\geq0$, (\ref{eq:lt}) still holds.
\end{lemma}
\begin{proof}
 Let us recall that the infinitesimal generator $L^0$ of the  process $X$ killed at time $T_0$ is given according to Theorem 6.1 in Lamperti \cite{Lamperti-72} by
  \begin{equation}\label{eq:inf_kill}
  L^0 f(x)= x^{-\alpha} L^{\xi} (f\circ \exp)(\log(x)),\quad x>0,
  \end{equation}
for at least functions $f$ such that $f(x),xf'(x),x^2f''(x)$ are continuous on $(0,\infty)$.  Next we have on the one the hand, writing $P^0_t$ for the semigroup  associated to $L^0$ and $\mathbf{I}(x)=1,x>0$,
\begin{eqnarray*}
L^0 \mathbf{I}(x) &=& \lim_{t\downarrow 0} \frac{P^0_t\mathbf{I}(x) -1}{t}\\
&=& \lim_{t\downarrow 0} \frac{\Q_x(T_0> t) -1}{t}\\
&=& -\lim_{t\downarrow 0} \frac{\Q_x(T_0\leq t)}{t}.
\end{eqnarray*}
On the other hand, we easily see,  from  \eqref{eq:inf_kill}, that
\[ L^0 \mathbf{I}(x)= x^{-\alpha}\psi(0),\quad x>0.\]
The first claim now follows.

Next, denote by $U^q$ (resp. $U^q_0$) the resolvent operator associated to the recurrent extension $X$ (resp. the process $X$ killed at time $T_0$). Then, an application of the strong Markov property yields the following identity 
\begin{equation}
U^{q}f(x)=U^{q}_0 f(x)+{\rm{E}}_x[e^{-qT_0}] U^qf(0). 
\end{equation}
Now note the following limits. Firstly, from the definition of the resolvent and the Feller property of $X$, $\lim_{q\rightarrow \infty } q U^qf(0)= f(0)$. Secondly, from
classical semi-group theory, see e.g.~\cite[Lemma 3.3]{Pazy-83},
\[\lim_{q\rightarrow \infty }q^2 U^{q}f(x)-qf(x)= L^{\psi} f(x) \textrm{ and } \lim_{q\rightarrow \infty }q^2 U_0^{q}f(x)-qf(x)= L^{0} f(x).
\]
Finally from the classical Tauberian Theorem,
\[ \lim_{q\rightarrow \infty }q {\rm{E}}_x[e^{-qT_0}]= \lim_{t\downarrow 0} \frac{\Q_x(T_0\leq t)}{t}. 
\]
The analytical expression of $L^{\psi} f(x)$ for any $x>0$ now follows. The boundary condition $\lim_{x\downarrow 0}x^{\theta-1}f'(x)=0$ is obtained by following a line of reasoning similar to Proposition 1.1 in Patie \cite{Patie-06c}. Moreover, the expression of the Laplace transform of $T_a$ is also deduced from arguments similar to Theorem 2.1 in \cite{Patie-06c}.
\end{proof}

Before proceeding  with the proof, let us make some remarks and  introduce some more notation.
Note that, by Lemma \ref{lem:transf}, the assumption $\alpha>0$  ensures that $\mathcal{T}_{\alpha}\psi$ is the Laplace exponent of a spectrally negative L\'evy process without killing.  Moreover, since $\theta<\alpha$ and $\psi$ is strictly convex, it follows that $\psi(\alpha)$  is strictly positive. Hence  $(\mathcal{T}_\alpha\psi)'(0^+) = \psi(\alpha)/\alpha>0$ which implies that the L\'evy process corresponding to $\mathcal{T}_\alpha\psi$ drifts to $+\infty$. Moreover, this also   implies that ${\rm P}^{\mathcal{T}_\alpha \psi}_0$ is necessarily the law of a transient positive self-similar Markov process with an entrance law at $0$.

 As we shall be dealing with first passage problems for $X$ (and hence also for $\xi$) we will make use of the so-called scale function $W_{\mathcal{T}_\alpha\psi}$ which satisfies $W_{\mathcal{T}_\alpha\psi}(x) =0$ for $x<0$ and otherwise is defined as the unique continuous function on $[0,\infty)$ with the Laplace transform
\[
\int_0^\infty e^{-u x}W_{\mathcal{T}_\alpha\psi}(x)d x = \frac{1}{\mathcal{T}_\alpha\psi(u)}\text{  for } u>0.
\]
Chapter 8 of Kyprianou \cite{Kyprianou-06}, Kyprianou and Palmowski \cite{Kyprianou-Palmowski-05}  and Chan et al. \cite{CKS} all expose analytical properties of $W_{\mathcal{T}_\alpha\psi}$  as well as  many fluctuation identities in which the scale function appears; some of which will be used below without further reference.
The proof of Theorem \ref{thm:ct} requires the following preliminary which concerns the quantity
\begin{eqnarray*}
\Oc{x}{a}= {\rm{E}}_x^{\mathcal{T}_\alpha\psi} \left[e^{-q \int_0^{\infty} \Id{X_s\leq
 a}ds}\right]
 \end{eqnarray*}
defined for any $x,q \geq 0$, $a>0$ and any Laplace exponent of a (possibly killed) spectrally negative L\'evy process, $\psi$.

\begin{lemma} \label{IDEsoln} Fix $q\geq 0$ and $\alpha>0$. Under the assumptions of Theorem \ref{thm:ct}
we have for all $x\geq 0$ and $a>0$
\begin{eqnarray}
\Oc{x}{a}&=&\frac{\mathcal{I}_{\mathcal{T}_\alpha\psi, \alpha} (q   (x\wedge a)^\alpha )}{\mathcal{I}_{\psi, \alpha} (q a^\alpha )}
- \frac{q a^\alpha }{\mathcal{I}_{\psi, \alpha} (q a^\alpha )}  \int_1^{\frac{x}{a}\vee 1} z^{-\alpha-1} W_{\mathcal{T}_\alpha\psi}(\log(z)) \mathcal{I}_{\mathcal{T}_\alpha\psi,\alpha}(q a^\alpha  z^{-\alpha}) dz\notag\\
&& +  \frac{\psi(\alpha)}{\alpha}W_{\mathcal{T}_\alpha\psi}(\log( x/a\vee 1))\left( 1- \frac{\mathcal{I}_{\psi, \alpha} (q a^{2\alpha}  (x\vee a)^{-\alpha})}{\mathcal{I}_{\psi, \alpha} (q a^\alpha )}\right).
\label{main}
\end{eqnarray}

\end{lemma}

 \begin{proof} 
From the self-similarity of $X$, we observe that the following
identity
\begin{equation}
\Oc{x}{a}=O^{\mathcal{T}_\alpha\psi}_{qa^{\alpha}}(x/a;1)
\label{eq:identity}
\end{equation}
is valid for any $a>0$ and $x\geq0$.
It therefore suffices to prove the identity for $a=1$.

We start by computing $\Oc{1}{1}$.
Let
\begin{eqnarray*}
\tau_1 = \inf \{s>0;\: X_s<1\}
 \end{eqnarray*}
 be the first passage time of $X$ below the level $1$.
 Fixing $y>1$ we may make use of the strong Markov property and spectral negativity to deduce that
\begin{align}
\Oc{1}{1} &= {\rm{E}}_1^{\mathcal{T}_\alpha\psi} \left[e^{-q \int_0^{T_y} \Id{X_s\leq
 1}ds}\right] {\rm{E}}_y^{\mathcal{T}_\alpha\psi} \left[e^{-q \int_0^{\infty} \Id{X_s\leq
 1}ds}\right] \notag\\
 &={\rm{E}}_1^{\mathcal{T}_\alpha\psi} \left[e^{-q \int_0^{T_y} \Id{X_s\leq
 1}ds}\right]
 \left( {\rm{E}}_y^{\mathcal{T}_\alpha\psi}\left[\Id{\tau_1<\infty} {\rm{E}}^{\mathcal{T}_\alpha\psi}_{X_{\tau_1}}\left[e^{-q T_1}\right]\right] \Oc{1}{1}+ {\rm{E}}_y^{\mathcal{T}_\alpha\psi}\left[\Id{\tau_1=\infty}\right]\right).
 \label{pathdecomp}
\end{align}
Solving for $ \Oc{1}{1}$ we get
\begin{equation}
 \Oc{1}{1}  =
 \frac{
 {\rm{E}}_y^{\mathcal{T}_\alpha\psi}\left[\Id{\tau_1=\infty}\right]
  }{
\left\{{\rm{E}}_1^{\mathcal{T}_\alpha\psi} \left[e^{-q \int_0^{T_y} \Id{X_s\leq
 1}ds}\right]\right\}^{-1}
 -   {\rm{E}}_y^{\mathcal{T}_\alpha\psi}\left[\Id{\tau_1<\infty} {\rm{E}}^{\mathcal{T}_\alpha\psi}_{X_{\tau_1}}\left[e^{-q T_1}\right]\right]
 }.
 \label{master}
\end{equation}

Now we evaluate some of the expressions on the right hand side above. First, we write  $\tau^\xi_0 = \inf\{s>0; \: \xi_s <0\}$.  On the one hand, recalling that $({\mathcal{T}_\alpha\psi})'(0^+) = \psi(\alpha)/\alpha>0$, we observe that
\begin{eqnarray}
 {\rm{E}}_y^{\mathcal{T}_\alpha\psi}\left[\Id{\tau_1=\infty}\right] &=&
 {\P}_{\log y}^{\mathcal{T}_\alpha\psi}\left(\int_0^{\tau^{\xi}_0}e^{\alpha\xi_s}ds=\infty\right)\notag\\
 &=&
 {\P}_{\log y}^{\mathcal{T}_\alpha\psi}\left(\tau^{\xi}_0 =\infty\right)\notag\\
 &=&\frac{\psi(\alpha)}{\alpha}W_{\mathcal{T}_\alpha\psi}(\log y)
 \label{incorporate too}
 \end{eqnarray}
 where the last line follows from the classical identity for the ruin probability in terms of scale functions (see for example Theorem 8.1 of \cite{Kyprianou-06}).
On the other hand, by Fubini's theorem (recalling the positivity of both the coefficients in the definition of $\mathcal{I}_{\mathcal{T}_\alpha\psi,\alpha}$ and of $X$), we have with the help of (\ref{eq:lt}),
\begin{eqnarray}
{\rm{E}}_y^{\mathcal{T}_\alpha\psi}\left[  \Id{\tau_1<\infty}{\rm{E}}^{\mathcal{T}_\alpha\psi}_{X_{\tau_1}}\left[e^{-q
T_1}\right]\right]&=&
{\rm{E}}_y^{\mathcal{T}_\alpha\psi}\left[\frac{\mathcal{I}_{\mathcal{T}_\alpha\psi,\alpha}(qX_{\tau_1}^{\alpha})\Id{\tau_1<\infty}}{\mathcal{I}_{\mathcal{T}_\alpha\psi,\alpha}(q)}\right]\notag\\
&=& \frac{1}{\mathcal{I}_{\mathcal{T}_\alpha\psi,\alpha}(q)}
{\rm{E}}_y^{\mathcal{T}_\alpha\psi}\left[\sum_{n=0}^{\infty}a_n(\mathcal{T}_\alpha\psi ; \alpha)
q^nX_{\tau_1}^{\alpha
n}\Id{\tau_1<\infty}\right]\notag\\
&=& \frac{1}{\mathcal{I}_{\mathcal{T}_\alpha\psi,\alpha}(q)} \sum_{n=0}^{\infty}a_n(\mathcal{T}_\alpha\psi ; \alpha) q^n
\E_{\log y}^{\mathcal{T}_\alpha\psi}\left[e^{\alpha n\xi_{\tau^{\xi}_0}}\Id{\tau^{\xi}_0<\infty}\right].
\label{sum}
 \end{eqnarray}

Next  we recall a known identity for spectrally negative L\'evy processes. Namely that for $x\geq 0$ and $u\geq0$, taking account of the fact that $({\mathcal{T}_\alpha\psi})'(0^+)> 0$,
\begin{equation}
\E^{\mathcal{T}_\alpha\psi}_{x}(e^{u \xi_{\tau^\xi_0}} \Id{\tau^{\xi}_0<\infty} ) = e^{ux} -\mathcal{T}_\alpha\psi(u)e^{ux}\int_0^x e^{-u z}W_{\mathcal{T}_\alpha\psi}(z)dz -\frac{\mathcal{T}_\alpha\psi(u)}{u}W_{\mathcal{T}_\alpha\psi}(x),
\label{incorporate}
\end{equation}
where $\mathcal{T}_\alpha\psi(u)/u$ is understood to be $(\mathcal{T}_\alpha\psi)'(0^+)$ when $u=0$.
See e.g.~\cite{Kyprianou-Palmowski-05}.
Hence incorporating (\ref{incorporate too}), (\ref{sum})
and (\ref{incorporate}) into (\ref{master}),
 recalling the identity (\ref{eq:lt}), and then taking limits as $y\downarrow 1$,  we have
\begin{align*}
 \lefteqn{\frac{1}{\Oc{1}{1}}}
& \\
&=\lim_{y\downarrow 1}
\frac{\alpha}{\psi(\alpha) \mathcal{I}_{\mathcal{T}_\alpha\psi,\alpha}(q)}
\sum_{n=0}^{\infty} a_n(\mathcal{T}_\alpha\psi ; \alpha)q^n
\left\{  \mathcal{T}_\alpha\psi(\alpha n) y^{\alpha n}\int_0^{\log y} e^{-\alpha n z}\frac{W_{\mathcal{T}_\alpha\psi}(z)}{W_{\mathcal{T}_\alpha\psi}(\log y)}dz + \frac{\mathcal{T}_\alpha\psi(\alpha n)}{\alpha n}\right\}\\
&-\frac{\alpha \mathcal{I}_{\mathcal{T}_\alpha\psi,\alpha}(qy^\alpha)}{\psi(\alpha) \mathcal{I}_{\mathcal{T}_\alpha\psi,\alpha}(q)}
{\rm{E}}_1^{\mathcal{T}_\alpha\psi} \left[e^{-q \int_0^{T_y} \Id{X_s\leq 1}ds}\right]^{-1}
\frac{1}{W_{\mathcal{T}_\alpha\psi}(\log y)}
\left\{ {\rm{E}}_1^{\mathcal{T}_\alpha\psi} \left[e^{-q \int_0^{T_y} \Id{X_s\leq 1}ds}\right]- {\rm E}_1^{\mathcal{T}_\alpha\psi}[e^{-q T_y}] \right\}.
\end{align*}
Now using the simple estimate, for $\theta\geq \epsilon\geq 0$,  $e^{-(\theta -\epsilon)}-e^{-\theta}\leq \epsilon$, we have that
\begin{eqnarray*}
{\rm{E}}_1^{\mathcal{T}_\alpha\psi} \left[e^{-q \int_0^{T_y} \Id{X_s\leq 1}ds}\right]- {\rm E}_1^{\mathcal{T}_\alpha\psi}[e^{-q T_y}]&\leq & q  {\rm E}_1^{\mathcal{T}_\alpha\psi}\left[ \int_0^{T_y} \Id{X_s>1}ds \right]\\
&=&  q  \mathbb{E}_0^{\mathcal{T}_\alpha\psi}\left[\int_0^{T^\xi_{\log  y}}  e^{\alpha\xi_s} \Id{\xi_s >0}ds\right]\\
&=&   q  \int_0^{\log y} e^{\alpha z} u(\log y, z)dz\\
&\leq&   q   y^\alpha \int_0^{\log y} u(\log y, z)dz,
\end{eqnarray*}
where $T^\xi_{\log y} =\inf\{s>0; \: \xi_s = \log y\}$ and $u(\log y, z)$ is the potential density of $(\xi, \mathbb{P}_0^{\mathcal{T}_\alpha\psi})$ when killed on exiting $(-\infty,\log y)$.  It can easily be deduced from Theorem 8.7, and the discussion that follows it, in \cite{Kyprianou-06} that, for $z>0$, $u(\log y, z) = W_{\mathcal{T}_\alpha\psi}(\log y - z)$. 
Now
note that $W_{\mathcal{T}_\alpha\psi}$ is monotone increasing and hence
\[
\lim_{y\downarrow 1}\frac{1}{W_{\mathcal{T}_\alpha\psi}(\log y)}
\left\{ {\rm{E}}_1^{\mathcal{T}_\alpha\psi} \left[e^{-q \int_0^{T_y} \Id{X_s\leq 1}ds}\right]- {\rm E}_1^{\mathcal{T}_\alpha\psi}[e^{-q T_y}] \right\}\leq \lim_{y\downarrow 1}  q   y^\alpha \int_0^{\log y} \frac{W_{\mathcal{T}_\alpha\psi}(\log y - z)}{W_{\mathcal{T}_\alpha\psi}(\log y)}dz  = 0.
\]
Similarly we have
\[
\lim_{y\downarrow 1}\sum_{n=0}^{\infty} a_n(\mathcal{T}_\alpha\psi ; \alpha)q^n
 \mathcal{T}_\alpha\psi(\alpha n) y^{\alpha n}\int_0^{\log y} e^{-\alpha n z}\frac{W_{\mathcal{T}_\alpha\psi}(z)}{W_{\mathcal{T}_\alpha\psi}(\log y)}dz \leq \lim_{y\downarrow 1} \mathcal{I}_{\mathcal{T}_\alpha\psi,\alpha}(qy^\alpha) \log y= 0.
\]
With the trivial observation that $\lim_{y\downarrow 1}{\rm{E}}_1^{\mathcal{T}_\alpha\psi} \left[e^{-q \int_0^{T_y} \Id{X_s\leq 1}ds}\right]=1,$ we thus conclude that
\[
 \Oc{1}{1}= \frac{\psi(\alpha) \mathcal{I}_{\mathcal{T}_\alpha\psi, \alpha}(q)}{ \alpha\sum_{n=0}^{\infty}
 a_n(\mathcal{T}_\alpha\psi ; \alpha)
 q^n
\frac{\mathcal{T}_\alpha\psi(\alpha n)}{\alpha n}
}.
\]

For any $n\geq1$,
\begin{eqnarray}
\label{p(x)/x}
\frac{\psi(\alpha) \alpha n}{\alpha\mathcal{T}_\alpha\psi(\alpha n)} a_n(\mathcal{T}_\alpha\psi ; \alpha)^{-1}
 =  \prod_{k=1}^n
\psi(\alpha k )
\end{eqnarray}
and with the interpretation that $\mathcal{T}_\alpha\psi(u)/u=(\mathcal{T}_\alpha\psi)'(0^+)$ when $u=0$, we also see that the left hand side of (\ref{p(x)/x}) is also equal to $1$ when $n=0$.
We finally come to rest at  the identity
\begin{equation}
 \Oc{1}{1} = \frac{ \mathcal{I}_{\mathcal{T}_\alpha\psi, \alpha}(q)}{\sum_{n=0}^{\infty}a_n(\psi;\alpha)
q^n
}
= \frac{ \mathcal{I}_{\mathcal{T}_\alpha\psi, \alpha}(q)  }{\mathcal{I}_{\psi, \alpha} (q)}.
\label{O(1)}
\end{equation}

Next, note that for all $0\leq x\leq 1$,
\[
 \Oc{x}{1}  = {\rm E}_x^{\mathcal{T}_\alpha \psi}(e^{-q T_1})\Oc{1}{1}
\]
and therefore, taking account of (\ref{eq:lt}) and Lemma \ref{lem:lt} we get
the expression given in
 \eqref{main} when $x\leq 1$ and $a=1$.

To get an expression when  $x>1$, we proceed as in (\ref{pathdecomp}) and we note from  (\ref{incorporate too}),  (\ref{p(x)/x}), (\ref{O(1)}) and Fubini's Theorem that
\begin{align}
{\rm E}^{\mathcal{T}_\alpha\psi}_x\left(e^{-q\int_0^\infty\Id{X_s \leq 1} ds}\right)
&\notag\\
&\hspace{-3cm}=  {\rm{E}}_x^{\mathcal{T}_\alpha\psi}\left[\Id{\tau_1<\infty} {\rm{E}}^{\mathcal{T}_\alpha\psi}_{X_{\tau_1}}\left[e^{-q T_1}\right]\right] \Oc{1}{1}+ {\rm{E}}_x^{\mathcal{T}_\alpha\psi}\left[\Id{\tau_1=\infty}\right]\notag\\
&\hspace{-3cm}=   \frac{1}{\mathcal{I}_{\psi, \alpha} (q)} \sum_{n=0}^{\infty} a_n(\mathcal{T}_\alpha\psi; \alpha) q^n
\left[1 -\mathcal{T}_\alpha\psi(\alpha n)\int_0^{\log x} e^{-\alpha n z}W_{\mathcal{T}_\alpha\psi}(z)dz -\frac{\mathcal{T}_\alpha\psi(\alpha n)}{\alpha n}x^{-\alpha n}W_{\mathcal{T}_\alpha\psi}(\log x)
\right]\notag\\
&+\frac{\psi(\alpha)}{\alpha}W_{\mathcal{T}_\alpha\psi}(\log x)\notag\\
&\hspace{-3cm}=  \frac{\mathcal{I}_{\mathcal{T}_\alpha\psi, \alpha}(q)}{\mathcal{I}_{\psi, \alpha} (q)}  - \frac{q}{\mathcal{I}_{\psi, \alpha} (q)}  \int_0^{\log x} e^{-\alpha z}W_{\mathcal{T}_\alpha\psi}(z) \mathcal{I}_{\mathcal{T}_\alpha\psi,\alpha}(q e^{-\alpha z}) dz\notag\\
&+\frac{\psi(\alpha)}{\alpha}W_{\mathcal{T}_\alpha\psi}(\log x)\left( 1- \frac{\mathcal{I}_{\psi, \alpha} (q x^{-\alpha})}{\mathcal{I}_{\psi, \alpha} (q)}\right)
\label{x>1}
\end{align}
which, after a change of variable, also agrees with
 (\ref{main}) when $a=1$.
 \end{proof}

The proof of our main theorem is now a very straightforward argument.  Indeed, note that the Laplace transform of $\Oc{0}{a}$  coincides with the one of the stopping time $(T_a, \Q^{\psi}_0)$  as given in (\ref{eq:lt}) and Lemma \ref{lem:lt}. The claim of equality in distribution follows from the injectivity of the Laplace transform. Finally, the self-decomposability of the pair follows from the proved self-decomposability of  $(T_a, \Q^{\mathcal{T}_{\alpha}\psi}_0)$  in Theorem 2.6 of \cite{Patie-06c}.

\section{Examples}

We refer to Lebedev's monograph
\cite{Lebedev-72} for detailed information on the special functions appearing in the examples below.

\begin{exa}[The case of no jumps: Bessel processes]
 When $\psi$ has no jump component it is possible to extract the original Ciesielski-Taylor identity for Bessel processes from Theorem \ref{thm:ct}. Indeed, we may take   $\alpha = 2$ and
\[
\psi_\nu(u)= \frac{1}{2}u^2 + \left(\frac{\nu}{2}-1\right) u
\]
where $\nu>0$. In that case it follows that ${\rm P}^\psi_\cdot$ is the law of a Bessel process of dimension $\nu$ as described in the introduction. Note that the root $\theta$ is zero for $\nu\geq 2$ and when $\nu\in(0,2)$ we have $\theta = 2-\nu<2$ thereby fulfilling the conditions of Theorem \ref{thm:ct}. The transformation $\mathcal{T}_{2}$ gives us the new Laplace exponent
\[
\mathcal{T}_{2}\psi_\nu (u) = \frac{1}{2}u^2 + \frac{\nu}{2} u  = \psi_{\nu+2}(u).
\]
The identity (\ref{eq:identity}) therefore agrees with the original identity (\ref{eq:ict}).
It is straightforward to show that
 \begin{eqnarray*}
 \I_{\psi_{\nu}, 2}(x) &=&  \Gamma(\nu/2) I_{\nu/2 -1} (\sqrt{2x})(\sqrt{2x} /2)^{-(\nu/2-1) }
 \end{eqnarray*}
 where
\[
{\rm{I}}_{\gamma}(x)=\sum_{n=0}^{\infty}\frac{(x/2)^{\gamma+2n}}{n!\Gamma(\gamma+n+1)}
\]
 stands for the modified Bessel function of index $\gamma$.
Hence  it follows that the shared Laplace transform on both sides of the identity is given by
both left and right hand side of these identities have Laplace transform given by
\[ \frac{1}{\I_{\psi_{\nu}, 2}(qa^2)}=\frac{(a\sqrt{2q})^{\nu/2-1}}{2^{\nu/2-1}\Gamma(\nu/2){\rm{I}}_{\nu/2-1}\left(a\sqrt{2q}\right)},
\]
thereby agreeing with the Laplace transform for the classical Ciesielski-Taylor  identity for Bessel processes.
Let us again fix $\alpha=2$ and   consider now, for $\kappa>0$, the Laplace exponent
\begin{eqnarray*}
\psi_{\nu,\kappa}(u)&=& \psi_{\nu}(u) -\kappa\\
&=& \frac{1}{2}(u-\theta_+)(u-\theta_-)
\end{eqnarray*}
where $\theta_{\pm}=1-\frac{\nu}{2}\pm \sqrt{(\frac{\nu}{2}-1)^2+2\kappa}$. We easily verify that under the additional condition $\kappa<\nu$, we have $\theta_+<2$ and thus ${\rm P}^\psi_\cdot$ stands for the law of the recurrent extension of a Bessel process of dimension $\nu$ killed at a rate $\kappa A_t$; recall that $A_t$ was defined in (\ref{A_t}). From the linearity property of the mapping $\mathcal{T}$, we get that
\[
\mathcal{T}_{2}\psi_{\nu,\kappa}(u) =  \psi_{\nu+2}(u)-\kappa\frac{u}{u+2}=\frac{u}{2(u+2)}(u+2-\theta_+)(u+2-\theta_-)
\]
which is the Laplace exponent of a linear Brownian motion with independent compound Poisson, exponentially distributed, negative jumps. The Ciesielski-Taylor identity now has the interesting feature that the process $(X, {\rm P}_0^{\psi_{\nu, \kappa}})$ has no jumps in its path prior to its moment of  reaching $0$, however the process $(X, {\rm P}_0^{\mathcal{T}_2\psi_{\nu,\kappa}})$ experiences  discontinuities with finite activity.

In this case we may also compute
\begin{eqnarray*}
 \I_{\psi_{\nu,\kappa},2}(x) &=&  {}_{1}F_2(1;1-\theta_+/2,1-\theta_-/2;x/2)\\
 \I_{\mathcal{T}_{2}\psi_{\nu,\kappa},2}(x) &=&  {}_{1}F_2(2;2-\theta_+/2,2-\theta_-/2;x/2)
 \end{eqnarray*}
 for the hypergeometric function ${}_{1}F_2(\delta;\beta,\gamma;x)=\sum_{n=0}^{\infty}\frac{(\delta)_{n}}{(\beta)_{n}(\gamma)_{n}n!}x^n$, with $(\delta)_{n}=\frac{\Gamma(\delta+n)}{\Gamma(\delta)}$.
 \end{exa}

\begin{exa}[The spectrally negative $\mathcal{T}$-Lamperti stable process]\label{TSexcluded}
 Caballero and Chaumont \cite{Caballero-Chaumont}  determined the characteristic triplet of the underlying L\'evy processes associated via the Lamperti mapping to $\alpha$-stable L\'evy  processes killed upon entering the negative half-line as well as two different  $h$-transforms thereof. In \cite{Patie-CBI-09}  and \cite{Chaumont-Kyprianou-Pardo-09} the Laplace exponent of these L\'evy processes in the spectrally negative case have been computed. In particular, recalling the notation $(u)_\alpha = \Gamma(u+\alpha)/\Gamma(u),$  the L\'evy process underlying  the $\alpha$-stable process killed upon entering into $(-\infty,0)$, is determined by the following Laplace exponent, for any $ 1<\alpha< 2$,
\begin{equation} \label{eq:lap_poch}
 \psi_{\alpha}(u)=c (
u+1-\alpha)_{\alpha}, \quad u\geq0,
\end{equation}
where $c$ is a positive constant which, for sake of simplicity, we set to $1$.
Note that $\psi_{\alpha}(0)=(1-\alpha)_{\alpha}<0$, showing that $\psi_\alpha$ corresponds to a killed spectrally negative L\'evy process and $\psi_{\alpha}(\alpha-1)=0$, showing that $\theta<\alpha$. Thus, the conditions of Theorem \ref{thm:ct} being satisfied, we recover the Ciesielski-Taylor identity
with  $\mathcal{T}_{\alpha}\psi_\alpha(u) = (
u)_{\alpha}$. Taking account of the fact that the latter Laplace exponent is that of the spectrally negative L\'evy process found in the Lamperti transformation which describes a spectrally negative $\alpha$-stable process conditioned to stay positive,
the Ciesielski-Taylor identity may otherwise be read as saying the following. For $1<\alpha<2$,  the law of the first passage time to unity of a spectrally negative $\alpha$-stable process reflected in its infimum
is equal to that of the occupation of the same process conditioned to stay positive. When seen in this context, we also see that we have recovered the only self-similar case of the identity mentioned in the remark at the bottom of p.~1475 in Bertoin \cite{BertoinCT}.

We also note for this example that
\begin{eqnarray*}
\mathcal{I}_{\psi_\alpha, \alpha}(x)&=&E_{\alpha,1}\left(
x\right)\\
\mathcal{I}_{\mathcal{T}_{\alpha}\psi_\alpha, \alpha}(x)&=&E_{\alpha,\alpha}\left(
x\right)
\end{eqnarray*}
 where $E_{\alpha,\beta}(x)=\sum_{n=0}^{\infty}\frac{1}{(\beta)_{\alpha n}}x^n$ stands for the Mittag-Leffler function.
 More generally, we introduce, for any $\beta >0$, the two-parameters family of L\'evy processes having the Laplace exponent $\psi_{\alpha,\beta}=\mathcal{T}_{\beta} \psi_{\alpha}$, that is
\[\psi_{\alpha,\beta}(u)=\frac{u}{u+\beta}(
u+\beta+1-\alpha)_{\alpha}.\]
We refer to them as $\mathcal{T}$-Lamperti stable processes. We easily check that $\psi_{\alpha, \beta}'(0^+)=(
\beta+1-\alpha)_{\alpha}\geq 0$ if $\beta\geq \alpha-1$ and otherwise $\psi_{\alpha, \beta}(\alpha-\beta-1)=0$. Thus, for any $\beta>0$, we get, from the composition property of the transformation $\mathcal{T}$, the identity in distribution
 \begin{equation}
\left(T_a, {{\rm P}}^{\psi_{\alpha,\beta}}_0\right)
\stackrel{(d)}{=} \left(\int_0^{\infty}\Id{X_s\leq a} ds,
{{\rm P}}^{\psi_{\alpha,\beta+\alpha}}_0\right).
\end{equation}
  Moreover, it is straightforward to verify that
\begin{equation*}
\mathcal{I}_{\psi_{\alpha,\beta}, \alpha}(x)={}_1\Psi_1\left( \left.\begin{array}{c}
                  (1,\beta/\alpha) \nonumber \\
                  \left(\alpha,\beta\right)
                \end{array} \right|
x\right)
\end{equation*}
where ${}_1\Psi_1\left( \left.\begin{array}{c}
                  (a,\alpha) \nonumber \\
                  \left(b,\beta\right)
                \end{array} \right|
x\right)=\sum_{n=0}^{\infty}\frac{(\alpha)_{an}}{(\beta)_{bn}n!}x^n$ stands for the Wright hypergeometric function.
\end{exa}


\begin{exa}[The spectrally negative saw-tooth process] \label{ex:st}
We consider the so-called saw-tooth process introduced and studied by Carmona et al.~\cite{Carmona-Petit-Yor-98}, and we would also like to note, the inspiration for the current article. It is a positive self-similar Markov process with index $1$  where the associated L\'evy process is the negative of the compound Poisson process of parameter
 $\kappa>0$ whose jumps are distributed as exponentials of parameter $\gamma+\kappa-2>0$ and a positive drift  of parameter  $1$, i.e.~\[\psi(u) =u\frac{u+\gamma -2 }{u+\gamma + \kappa -2}, \: u\geq0.\] Note that $\psi(0) = 0$ and  $\psi'(0^+)=\frac{\gamma -2 }{\gamma + \kappa -2}$ showing transience for  $\gamma\geq 2$ and recurrence for $\gamma\in(2-\kappa, 2)$ of the law ${\rm P}^\psi_0$. We also note that $\theta = 0$ if $\gamma\geq 2$ and otherwise $\theta \in(0, \kappa)$ when $\gamma\in(2-\kappa, 2)$. Therefore, we assume that either $\gamma\geq2$ and $\kappa>0$ or $\gamma\in(2-\kappa, 2)$ with $\kappa<1$.  Moreover,
\begin{equation*}
\mathcal{I}_{\psi,1}(x)={}_1F_1(\gamma+\kappa-1,\gamma-1;x)
\end{equation*}
where ${}_1F_1(a,b;x)=\sum_{n=0}^{\infty} \frac{(a)_n}{(b)_n}\frac{x^n}{n!}$ stands for the confluent hypergeometric function. Next note that
\begin{equation*}
\mathcal{T}_{1}\psi(u)= u\frac{u+\gamma -1 }{u+\gamma + \kappa -1},
\end{equation*}
and we get that $\mathcal{T}_{1}\psi$ is the Laplace exponent of the negative of the compound Poisson process of parameter
 $\kappa>0$ whose jumps are distributed as exponentials of parameter $\gamma+\kappa-1$ and a positive drift of  parameter $1$.
Finally, for any $a>0$, we have, appealing to  obvious notation,
\begin{equation}
\left(T_a, \Q^{\gamma-1,\kappa}_0\right)
\stackrel{(d)}{=} \left(\int_0^{\infty}\Id{X_s\leq a} ds,
\Q^{^{\gamma,\kappa}}_0\right).
\end{equation}
Note that, in terms of our notation, Carmona et al.~\cite[Theorem 4.8]{Carmona-Petit-Yor-98} obtain the identity
\begin{equation*}
\left(T_a, \Q^{\gamma-1,\kappa}_0\right)
\stackrel{(d)}{=} \left(\int_0^{\infty}\Id{X_s\leq a} ds,
\Q^{^{\gamma,\kappa-1}}_0\right)
\end{equation*}
which differs with the identity provided by our main result. After consulting with Carmona and Yor on this point, it appears that there is an error in their proof which explains the discrepancy. 

\bigskip

One important point which comes out of the analysis in Carmona et al.~\cite{Carmona-Petit-Yor-98} is the relation of the   solution obtained in Theorem \ref{thm:ct} with a certain integro-differential equation, even in the general setting of Theorem \ref{thm:ct}. Suppose that $L^{\mathcal{T}_\alpha\psi}$ is the infinitesimal generator of $(X, {\rm P}^{\mathcal{T}_\alpha\psi}_0)$.
Then a standard martingale argument  shows that if a non-negative solution to the integro-differential equation
\begin{eqnarray}\label{IDE}
L^{\mathcal{T}_\alpha\psi} u(x) = q\Id{x\leq 1} u(x)
\end{eqnarray}
exists for $q\geq 0$ (note that no boundary conditions are required at $0$ as ${\rm P}^{\mathcal{T}_\alpha\psi}_0$ is  an entrance law under the assumptions of Theorem \ref{thm:ct}) then
\[
u(x) = {\rm E}^{\mathcal{T}_\alpha\psi}_x\left(e^{-q\int_0^\infty\Id{X_s \leq 1} ds}\right) .
\]
The  issue of solving (\ref{IDE}) has been circumvented in our presentation by approaching the problem through fluctuation theory instead. One should therefore think of the expression  (\ref{IDEsoln})  as providing the solution to the integro-differential equation (\ref{IDE}). Note in particular that there is continuity in (\ref{IDEsoln}) at $x=1$; a requirement that was sought by Carmona et al. \cite{Carmona-Petit-Yor-98} when trying to solve the integro-differential equation explicitly in their setting.
\end{exa}


\section*{Acknowledgements} We would like to thank an anonymous referee for their helpful remarks on an earlier draft of this paper.


\begin{thebibliography}{99}

\bibitem{BertoinCT}
J.~Bertoin.
\newblock An extension of Pitman's theorem for spectrally positive L\'evy processes.
\newblock {\em Ann. Probab.} 20:1464--1483, 1992.

\bibitem{Bertoin-96}
J.~Bertoin.
\newblock {\em L\'evy Processes}.
\newblock Cambridge University Press, Cambridge, 1996.

\bibitem{Bertoin-Yor-02-b}
J.~Bertoin and M.~Yor.
\newblock The entrance laws of self-similar {M}arkov processes and exponential
  functionals of {L}\'evy processes.
\newblock {\em Potential Anal.}, 17(4):389--400, 2002.


\bibitem{Biane-85}
Ph. Biane.
\newblock Comparaison entre temps d'atteinte et temps de s\'ejour de certaines
  diffusions r\'eelles.
\newblock In {\em S\'eminaire de probabilit\'es, XIX, 1983/84}, volume 1123 of
  {\em Lecture Notes in Math.}, pages 291--296. Springer, Berlin, 1985.

\bibitem{Ca-Ch} M.E. ~Caballero and L. ~Chaumont.
\newblock Weak convergence of positive self-similar Markov processes and overshoots of L\'evy processes.
\newblock {\it Ann. Probab.} {\bf 34}, 1012--1034, 2006.


\bibitem{Caballero-Chaumont} M.E. ~Caballero and L. ~Chaumont.
\newblock Conditioned stable L\'evy processes and the Lamperti representation.
\newblock {\it J. Appl. Probab.},  {\bf 43}, 967--983, 2006.



\bibitem{Carmona-Petit-Yor-98}
Ph. Carmona, F.~Petit, and M.~Yor.
\newblock Beta-{gamma} random variables and intertwining relations between
  certain {Markov} processes.
\newblock {\em Rev. Mat. Iberoamericana}, 14(2):311--368, 1998.


\bibitem{CKS} T. Chan, A.E.~Kyprianou and M. Savov.
\newblock  Smoothness of scale functions for spectrally negative L\'evy processes
\newblock \textit{To appear in Probability Theory and Related Fields}, 2010.


\bibitem{Chaumont-Kyprianou-Pardo-09}
L. Chaumont, A.E. ~Kyprianou and J.C. ~Pardo.
\newblock Some explicit identities associated with positive self-similar Markov processes.
\newblock {\em Stoch. Proc. Appl.} {\bf 119}:980--1000, 2009.

\bibitem{Ciesielski-Taylor-62}
Z.~Ciesielski and S.J. Taylor.
\newblock First passage times and sojourn times for {B}rownian motion in space
  and the exact {H}ausdorff measure of the sample path.
\newblock {\em Trans. Amer. Math. Soc.}, 103:434--450, 1962.

\bibitem{Fitz} P.J. Fitzsimmons.
\newblock On the existence of recurrent extension of positive self-similar Markov processes.
\newblock {\it Elec. Comm. Probab.} 11: 230--241, 2006.

\bibitem{Getoor-Sharpe-79}
R.K. Getoor and M.J. Sharpe.
\newblock Excursions of {B}rownian motion and {B}essel processes.
\newblock {\em Z. Wahr.}, 47:83--106, 1979.

\bibitem{Hubalek-Kyprianou-10}
F. Hubalek and A.~E. Kyprianou.
\newblock Old and new examples of scale functions for spectrally negative L\'evy processes.
\newblock To appear in {\em Stochastic Analysis, Random Fields and Applications. Eds. R. Dalang, M. Dozzi and F. Russo}, Progress in Probability, Birkh\"auser, 2010.

\bibitem{Kyprianou-06}
A.~E.~Kyprianou.
\newblock {\em Introductory lectures on fluctuations of {L}\'evy processes with
  applications}.
\newblock Universitext. Springer-Verlag, Berlin, 2006.


\bibitem{Kyprianou-Palmowski-05}
A.E. Kyprianou and Z.~Palmowski.
\newblock A martingale review of some fluctuation theory for spectrally
  negative {L}\'evy processes.
\newblock In {\em S\'eminaire de Probabilit\'es {XXXVIII}}, volume 1857 of {\em
  Lecture Notes in Math.}, pages 16--29. Springer, Berlin, 2005.

\bibitem{Lamperti-72}
J.~Lamperti.
\newblock Semi-stable {M}arkov processes. {I}.
\newblock {\em Z. Wahrsch. Verw. Geb.}, 22:205--225, 1972.

\bibitem{Lebedev-72}
N.N. Lebedev.
\newblock {\em Special Functions and their Applications}.
\newblock Dover Publications, New York, 1972.


\bibitem{Patie-CBI-09}
P.~Patie.
\newblock Exponential functional of one-sided {L}\'evy processes and
  self-similar continuous state branching processes with immigration.
\newblock {\em Bull. Sci. Math.}, 133(4):355--382, 2009.

\bibitem{Patie-06c}
P.~Patie.
\newblock Infinite divisibility of solutions to some self-similar
  integro-differential equations and exponential functionals of {L}\'evy
  processes.
\newblock {\em Ann. Inst. H. Poincar\'e Probab. Statist.}, 45(3):667--684, 2009.


\bibitem{Pazy-83}
A.~Pazy.
\newblock {\em Semigroups of linear operators and applications to partial
  differential equations.}, volume~44 of {\em Applied Mathematical Sciences}.
\newblock Springer-Verlag, New York, Berlin, Heidelberg, Tokyo, 1983.


\bibitem{Rivero1}
V.~Rivero.
\newblock Recurrent extensions of self-similar {M}arkov processes and
  {C}ram\'er's condition.
\newblock {\em Bernoulli}, 11(3):471--509, 2005.


\bibitem{Yor-91}
M.~Yor.
\newblock Une explication du th\'eor\`eme de {Ciesielski-Taylor}.
\newblock {\em Ann. Inst. H. Poincar\'e Probab. Statist.}, 27(2):201--213, 1991.

\end{thebibliography}
\end{document}